\documentclass[12pt]{article}
\usepackage{graphicx}
\usepackage{subfig}
\usepackage{amsmath}
\usepackage{amsfonts}
\usepackage{amssymb}
\usepackage{amsthm}
\usepackage{newlfont}
\usepackage{color}
\usepackage{float}
\usepackage[font=small,labelfont=bf]{caption}
\begin{document}

\title{ Non-local Integrals and Derivatives on Fractal Sets with Applications}

\author{  Alireza K. Golmankhaneh $^{1\dag}$, Dumitru Baleanu  $^{2,3 \ddag}$ }

\maketitle \vspace{-9mm}

\begin{center}
$^1$ Young Researchers and Elite Club, Urmia Branch,\\ Islamic Azad university, Urmia, Iran.\\
\emph{$^\dag$E-mail address}: a.khalili@iaurmia.ac.ir\\

$^2$ Department of Mathematics \\
$\c{C}$ankaya University, 06530 Ankara, Turkey\\
$^3$Institute of Space
Sciences,\\
 P.O.BOX, MG-23, R 76900,~Magurele-Bucharest, Romania\\

\emph{$^\ddag$E-mail address}: dumitru@cankaya.edu.tr
\end{center}
\date{}
%%% ----------------------------------------------------------------------
\maketitle \vspace{-9mm}
%\begin{center}
%$^1$ Young Researchers and Elite Club, Urmia Branch, Islamic Azad University, Urmia, Iran,\\
% $^\dag$\emph{E-mail address}: a.khalili@iaurmia.ac.ir\\
%\end{center}
%%% ----------------------------------------------------------------------
\begin{abstract}
In this paper, we discussed the non-local derivative on the fractal Cantor set. The scaling properties are given for  both local and non-local fractal derivatives. The local and non-local fractal diﬀerential equations are solved and compared and related physical models are suggested.
\end{abstract}

 %\vspace{.5cm} \noindent {\it {Keywords:}}

\section{Introduction}
 Fractional calculus became an important tool which was applied successfully in many branches of science and engineering etc \cite{ab1,ab2,ab3,ab4,ab5}. The models based on fractional derivatives are crucial for describing the processes with memory effect \cite{ab6}. Local fractional has been defined on the real-line \cite{ab7}. As it is well known the integer, fractional and complex order derivatives and integrals are defined on the real-line. Analysis on the fractal has been studied by many researchers \cite{ab9,ab10,ab11}. The fractals curves and the functions on fractal space are not differentiable in the sense of  standard calculus. As a result, by this motivation  recently in the seminal paper the $F^{\alpha}$-calculus is suggested as a framework on the fractal sets and fractal curves \cite{ab12,ab13,ab14,ab15}. The $F^{\alpha}$-calculus is generalized and applied in physics as a new and useful tool for modelling   processes on the fractals. Newtonian mechanics and Schr$\ddot{o}$dinger equation on the fractal sets and curves are given \cite{ab16,ab17,ab18}. The gauge integral is utilized to generalized the $F^{\alpha}$-calculus for the unbound and singular function \cite{ab19}. The fractal grating is modeled by $F^{\alpha}$-calculus and corresponding diffraction is presented  \cite{ab19}. One of the important aspects of fractional calculus was transferred recently to the fractal derivatives. Namely, the concept of non-local fractal derivatives was introduced in \cite{Golmankhaneh-107-k}.~In this manuscript our main aim is to define the fractal non-local derivatives and study their properties.\\
The plane of this work is as follows:\\
In  Section \ref{2-sec} we summarize the basic definitions and properties of the the local fractional derivatives. In  Section \ref{3-sec} the scaling properties of local and non-local derivatives are derived. More, in Section \ref{4-sec} we develop the theory of fractal local and non-local  Laplace transformations. In Section  \ref{5-sec} the comparison of local and non-local linear fractal differential equations  are presented.~In Section \ref{6-sec} we indicate  some illustrative  applications. Section  \ref{7-sec} contains our conclusion.

\section{Preliminaries\label{2-sec}}
In this section we recall some basic definitions and properties of  the local fractal calculus (LFC) and non-local fractal calculus (NLFC) \cite{ab12,Golmankhaneh-107-k}.

\subsection{Local fractal calculus}
\begin{figure}[H]
  \centering
  \includegraphics[scale=0.5]{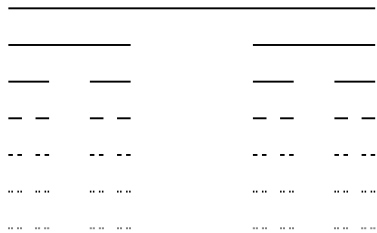}
  \caption{We present  triadic Cantor set by iteration.}\label{422vgtweqavdrb}
\end{figure}

\begin{figure}[H]
  \centering
  \includegraphics[scale=0.5]{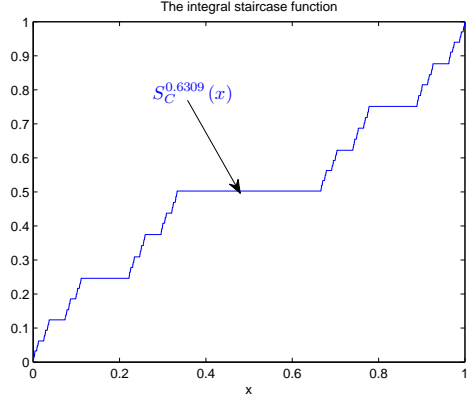}
  \caption{We indicate  the integral staircase function  for a triadic  Cantor set $F$. }\label{42des2vdrb}
\end{figure}

In the seminal paper local $F^{\alpha}$-calculus is built on fractal Cantor set which is shown in Figure [\ref{422vgtweqavdrb}] \cite{ab12}.
The integral staircase function $S_{F}^{\alpha}(x)$ of order
$\alpha$ for the triadic Cantor set $F$ is defined in  \cite{ab12} by
\begin{equation}\label{t}
    S_{F}^{\alpha}(x)=\begin{cases}
    \gamma^{\alpha}(F,a_{0},x) ~~~~\text{if} ~~~~~x\geq a_{0}\\
    -\gamma^{\alpha}(F,a_{0},x) ~~~~\text{otherwise},
\end{cases}
\end{equation}
where $a_{0}$ is an arbitrary real number. The graph of the integral staircase function is depicted in Figure [\ref{42des2vdrb}].
~Then $F^{\alpha}$-derivative is defined  for a function with this support as follows
 \cite{ab12}
\begin{equation}
    D_{F}^{\alpha}f(x)=\begin{cases}
    $F$-\lim_{y\rightarrow x} \frac{f(y)-f(x)}{S_{F}^{\alpha}(y)-S_{F}^
    {\alpha}(x)} ~~~\text{if} ~~~x\in F,\\
    0, ~~~~~~~~~~~~~\textmd{otherwise},
    \end{cases}
 \end{equation}
if the limit exists. For more details we refer the reader to \cite{ab12}.\\
\subsection{Non-local fractal calculus}
In this section, we review  the non-local derivatives and basic definitions \cite{Golmankhaneh-107-k}.\\
\textbf{Definition 1.} A function $f(S_{F}^{\alpha}(x)),~ x>0$ is in the space $\mathcal{C}_{F,\rho}, ~\rho\in \Re$ if there exists a real number $p>\rho$, such  $f(S_{F}^{\alpha}(x))=S_{F}^{\alpha}(x)^{p} f_{1}(S_{F}^{\alpha}(x))$, where $f_{1}(S_{F}^{\alpha}(x))\in \mathcal{C}_{F}^{\alpha}[a,b]$, and it is in the $\mathcal{C}_{F,\rho}^{n\alpha}[a,b]$ if and only if
\begin{equation}
(D_{F}^{\alpha})^{n} f(S_{F}^{\alpha}(x))\in \mathcal{C}_{F,\rho},~~~~ n\in N.
\end{equation}
Here and subsequently, we define the fractal left-sided  Riemann-Liouville  integral as follows
 \begin{align}\label{xsaz125}
 & _{a}\mathcal{I}_{x}^{\beta}f(x)\nonumber\\&:=\frac{1}{\Gamma^{\alpha}_{F}(\beta)}\int_{S_{F}^{\alpha}(a)}^ {S_{F}^{\alpha}(x)}\frac{f(t)}{(S_{F}^{\alpha}(x)-S_{F}^{\alpha}(t))^{\alpha-\beta}}d_{F}^{\alpha}t.
 \end{align}
 where $~S_{F}^{\alpha}(x)>S_{F}^{\alpha}(a)$.\\
\textbf{Definition 2.} The fractal left-sided Riemann-Liouville  derivative is defined as
\begin{align}\label{gcfrtgffvgf}
 & _{a}\mathcal{D}_{x}^{\beta} f(x)\nonumber\\&:=\frac{1}{\Gamma^{\alpha}_{F}(n-\beta)} (D_{F}^{\alpha})^n
\int_{S_{F}^{\alpha}(a)}^{S_{F}^{\alpha}(x)}\frac{f(t)} {(S_{F}^{\alpha}(x)-S_{F}^{\alpha}(t))^{-n\alpha+\beta+\alpha}}d_{F}^{\alpha}t.
 \end{align}
\textbf{Definition 3.} For A $f(x)\in C^{\alpha n} [a, b],~n\alpha-\alpha\leq\beta <\alpha n$ the fractal left-sided Caputo derivative is defined as
\begin{align}\label{gfvsewgf}
&  _{a}^{C}\mathcal{D}_{x}^{\beta} f(x)\nonumber\\&:=\frac{1}{\Gamma^{\alpha}_{F}(n-\beta)}
\int_{S_{F}^{\alpha}(a)}^{S_{F}^{\alpha}(x)}(S_{F}^{\alpha}(x)-S_{F} ^{\alpha}(t))^{n\alpha-\beta-\alpha} (D_{F}^{\alpha})^nf(t) d_{F}^{\alpha}t.
\end{align}
\textbf{Definition 4.} The fractal Gr\"{u}nwald and Marchaud derivative of a function $f(x)$ with support of fractal sets is defined as
\begin{align}\label{cdezaqw}
 & ^{G}\mathcal{D}^{\beta} f(x_{0})=\nonumber\\&F-\lim_{n \rightarrow \infty}\frac{1}{\Gamma_{F}^{\alpha}(-\beta)}\left(\frac{S_{F}^{\alpha}(x_{0})}{n}\right)^{-\beta}\sum_{k=0}^{n-1}
  \frac{\Gamma_{F}^{\alpha}(k-\beta)}{\Gamma_{F}^{\alpha}(k+1)}\nonumber\\&f\left(S_{F}^{\alpha}(x_{0})-k\frac{S_{F}^{\alpha}(x_{0})}{n}\right).\nonumber
\end{align}
\textbf{Definition 5.} The generalized fractal standard  Mittag-Leffler functions  is defined as \cite{Golmankhaneh-107-k}
\begin{equation}\label{dawsq9}
  E_{F,\eta}^{\alpha}(x)=\sum_{k=0}^{\infty}\frac{S_{F}^{\alpha}(x)^k}
  {\Gamma_{F}^{\alpha}(\eta k+1)},~~~\eta>0,~~\nu\in \Re.
\end{equation}
The fractal two parameter $\eta,~\nu$ Mittag-Liffler function is defined as
\begin{equation}\label{dawsq9}
  E_{F,\eta,\nu}^{\alpha}(x)=\sum_{k=0}^{\infty}\frac{S_{F}^{\alpha}(x)^k}
  {\Gamma_{F}^{\alpha}(\eta k+\nu)},~~~\eta>0,~~\nu\in \Re.
\end{equation}
\textbf{Definition 6.} For a given function $f(S_{F}^{\alpha}(x))$ the fractal Laplace transform is denoted by $F(s)$ and  defined  as \cite{Golmankhaneh-107-k}
\begin{equation}\label{sewa}
  \mathcal{F}_{F}^{\alpha}(S_{F}^{\alpha}(s))=\mathcal{L}_{F}^{\alpha}[f(x)]=
  \int_{S_{F}^{\alpha}(0)}^{S_{F}^{\alpha}(\infty)}f(x)
e^{-S_{F}^{\alpha}(s)S_{F}^{\alpha}(x)}d_{F}^{\alpha}x,
\end{equation}
where $S_{F}^{\alpha}(s)$ is limited by the values that the integral converges.~The function $f(S_{F}^{\alpha}(x))$ is $F$-continuous  and has following condition
\begin{equation}
\sup\frac{|f(S_{F}^{\alpha}(x))|}{e^{S_{F}^{\alpha}(c)S_{F}^{\alpha}(x)}}<\infty,~~~S_{F}^{\alpha}(c)\in \Re,~~~~S_{F}^{\alpha}(x)>0.
\end{equation}
In view of the above conditions the fractal Laplace transform exists for all $S_{F}^{\alpha}(s)> S_{F}^{\alpha}(c)$.  We follow the notation as $\mathcal{L}_{F}^{\alpha} [f(x)]= \mathcal{F}_{F}^{\alpha}(S_{F}^{\alpha}(s))$ and $\mathcal{L}_{F}^{\alpha} [g(x)]= \mathcal{G}_{F}^{\alpha}(S_{F}^{\alpha}(s))$.\\
\textbf{Remark 1.} We denote that if we choose $\beta=\alpha$ then we have
\begin{equation}\label{gfvgf}
  _{a}\mathcal{D}_{x}^{\alpha} f(x)= D_{F,x}^{\alpha}f(x)|_{x=S_{F}^{\alpha}(a)}.
 \end{equation}
\section{Scale properties of fractal local and non-local fractal calculus \label{3-sec} }
In this section we study  the scale properties of the LFC and NLFC.
\subsection{Scale change on the local fractal derivatives}
A function $f(S_{F}^{\alpha}(x))$ is called fractal homogenous of degree-$m\alpha$ or invariant under  fractal rescalings if we have
\begin{equation}\label{xdr}
  f(S_{F}^{\alpha}(\lambda x))=\lambda^{m \alpha}f(S_{F}^{\alpha}(x)),
\end{equation}
where for some $m$ and for all $\lambda$. The fractals have self-similar properties, namely for the case of function with the fractal  Cantor set support we choose $m=1$ and $\lambda=1/3^n, n=1,2,...$ then
\begin{equation}\label{xdrmn}
  f(S_{F}^{\alpha}(\frac{1}{3^n} x))=(\frac{1}{3^n})^{\alpha}f(S_{F}^{\alpha}(x)),
\end{equation}
where $\alpha=0.6 $ is the dimension of triadic Cantor set.~The fractal derivative of the fractal homogenous function  $f(S_{F}^{\alpha}(x))$   rescaling   as follows
\begin{equation}\label{xddcxsr}
  D_{F}^{\alpha}f(S_{F}^{\alpha}(\lambda x))=\lambda^{m\alpha-\alpha}f(S_{F}^{\alpha}(x)).
\end{equation}
\subsection{Scale change on the non-local fractal derivatives}
By a scale change of the fractal function $f(S_{F}^{\alpha}(x))$, we mean converts
\begin{equation}
x \rightarrow \lambda x \Rightarrow S_{F}^{\alpha}(\lambda x)=\lambda^{\alpha} S_{F}^{\alpha}(x),
\end{equation}
and using Eq. (\ref{gcfrtgffvgf}) and choosing $a=0$  we derive
\begin{equation}\label{xsaz258a}
  _{0}\mathcal{D}_{x}^{\beta}(f(S_{F}^{\alpha}(\lambda x)) )=\lambda^{\beta \alpha}~_{0}\mathcal{D}_{\lambda x}^{\beta}(f(S_{F}^{\alpha}(\lambda x)) ),
\end{equation}
which is called scale change on the non-local fractal derivatives.\\
\section{Laplace transformation on fractals \label{4-sec}}
Let us give some important lemmas  which are useful for finding the fractal Laplace transforms of function $f( S_{F}^{\alpha}(x))$.\\
\textbf{Lemma 1.} The fractal Laplace transform of the non-local  fractal Caputo derivative of order $m \alpha- \alpha<\beta\leq m\alpha$,~$m\in N$ is
\begin{align}\label{poi78}
 & \mathcal{L}_{F}^{\alpha}\{  _{0}^{C}\mathcal{D}_{x}^{\beta}f(x) \}=\frac{(S_{F}^{\alpha}(s))^{m \alpha}\mathcal{F}_{F}^{\alpha}(s)-(S_{F}^{\alpha}(s))^{m \alpha-\alpha}f(S_{F}^{\alpha}(0))}{S_{F}^{\alpha}(s)^{m\alpha-\beta}}\nonumber\\&
 \times \frac{-(S_{F}^{\alpha}(s))^{m \alpha-2 \alpha} D_{x}^{\alpha}f(x)|_{x=S_{F}^{\alpha}(0)}-\ldots - D_{x}^{m \alpha-\alpha}f(x)|_{x=S_{F}^{\alpha}(0)}}{1}.
\end{align}
\textbf{Proof:} We first compute the Laplace fractal transform of the fractal Caputo fractional derivative of order $\beta$ as follows
\begin{align}\label{bv}
  &\mathcal{L}_{F}^{\alpha}\{~_{0}^{C}\mathcal{D}_{x}^{\beta}f(x) \}\nonumber\\&=\mathcal{L}_{F}^{\alpha}\{_{0}\mathcal{I}_{x}^{m \alpha-\beta}  (D_{x}^{\alpha})^{m}f(x) \}\nonumber\\&=\frac{\mathcal{L}_{F}^{\alpha}[ (D_{x}^{\alpha})^{m}f(x)]}{s^{m\alpha-\beta}}
\end{align}
In view of Eq. (\ref{vcxzd}) which completes the proof.\\
\textbf{Lemma 2.} For a given $\zeta,~\mu>0,  S_{F}^{\alpha}(a)\in \Re$ and $S_{F}^{\alpha}(s)^{\zeta}> |S_{F}^{\alpha}(a)|$ the fractal Laplace transform is
\begin{align}\label{g895raws}
  & \mathcal{L}_{F}^{\alpha,-1}\left[\frac{S_{F}^{\alpha}(s)^{\zeta-\mu}}{S_{F}^{\alpha}(s)^{\zeta}+S_{F}^{\alpha}(a)}\right]\nonumber\\&=
   S_{F}^{\alpha}(x)^{\mu-1}E_{F,\zeta,\mu}^{\alpha}(-S_{F}^{\alpha}(a)S_{F}^{\alpha}(x)^{\zeta}).
 \end{align}
\textbf{Proof:} Using the series expansion we have
 \begin{eqnarray}\label{147}
 \frac{S_{F}^{\alpha}(s)^{\zeta-\mu}}{S_{F}^{\alpha}(s)^{\zeta}+S_{F}^{\alpha}(a)}&=&\frac{1}{S_{F}^{\alpha}(s)^{\mu}}
 \frac{1}{1+\frac{S_{F}^{\alpha}(a)}{S_{F}^{\alpha}(s)^{\zeta}}}\\ &=&\frac{1}{S_{F}^{\alpha}(s)^{\mu}}\sum_{n=0}^{\infty}
 \left(\frac{-S_{F}^{\alpha}(a)}{S_{F}^{\alpha}(s)^{\zeta}}\right)^n \nonumber\\
 &=&\sum_{n=0}^{\infty}
 \frac{(-S_{F}^{\alpha}(a))^n}{S_{F}^{\alpha}(s)^{n\zeta+\mu}}
 \end{eqnarray}
The inverse fractal Laplace transform of Eq. (\ref{147}) leads to
\begin{align}\label{x}
& \sum_{n=0}^{\infty}\frac{(-S_{F}^{\alpha}(a))^n~ S_{F}^{\alpha}(x)^{n\zeta+\mu-1}}{\Gamma_{F}^{\alpha}(n\zeta+\mu)} \nonumber\\&=S_{F}^{\alpha}(x)^{\mu-1}\sum_{n=0}^{\infty}\frac{(-S_{F}^{\alpha}(a)~S_{F}^{\alpha}(x)^{\zeta})^n}{\Gamma_{F}^{\alpha}(n\zeta+\mu)}\nonumber\\&
 =S_{F}^{\alpha}(x)^{\mu-1}E_{F,\zeta,\mu}^{\alpha}(-S_{F}^{\alpha}(a)S_{F}^{\alpha}(x)^{\zeta}).
\end{align}
\textbf{Lemma 3.} Suppose $\zeta\geq\mu>0$,~  $S_{F}^{\alpha}(a)\in \Re$ and $S_{F}^{\alpha}(s)^{\zeta-\mu}>|S_{F}^{\alpha}(a)|$ then we have
\begin{align}\label{nju159}
  &\mathcal{L}_{F}^{\alpha,-1}\left[\frac{1}{(S_{F}^{\alpha}(s)^{\zeta}+S_{F}^{\alpha}(a)S_{F}^{\alpha}(s)^{\mu})^{n+1}}\right]\nonumber\\&= S_{F}^{\alpha}(x)^{\zeta(n+1)-1}\sum_{k=0}^{\infty}\frac{-(S_{F}^{\alpha}(a))^k}{\Gamma_{F}^{\alpha}(k(\zeta-\mu)+(n+1)\zeta)} \nonumber\\& \binom{n+k}{k} S_{F}^{\alpha}(x)^{k(\zeta-\mu)}.
\end{align}
\textbf{Proof:} Let us use following expression
\begin{equation}\label{2}
  \frac{1}{(1+S_{F}^{\alpha}(x))^{n+1}}=\sum_{k=0}^{\infty} \binom{k+n}{k}(-S_{F}^{\alpha}(x))^k.
\end{equation}
Therefore we can write
\begin{align}
  &\frac{1}{(S_{F}^{\alpha}(s)^{\zeta}+S_{F}^{\alpha}(a)S_{F}^{\alpha}(s)^{\mu})^{n+1}} \nonumber\\&= \frac{1}{(S_{F}^{\alpha}(s)^{\zeta})^{n+1}}\frac{1}{(1+\frac{S_{F}^{\alpha}(a)}{S_{F}^{\alpha}(s)^{\zeta-\mu}})^{n+1}} \nonumber \\&=
   \frac{1}{(S_{F}^{\alpha}(s))^{n+1}}\sum_{k=0}^{\infty}\binom{n+k}{k}\left(\frac{-S_{F}^{\alpha}(a)}{S_{F}^{\alpha}(s)^{\zeta-\mu}}\right)^{k}.\nonumber
\end{align}
The proof is complete.\\
\textbf{Lemma 4.} For $\zeta\geq\mu,~\zeta>\xi,~S_{F}^{\alpha}(a)\in \Re,~S_{F}^{\alpha}(s)^{\zeta-\mu}>|S_{F}^{\alpha}(a)|$ and $|S_{F}^{\alpha}(s)^{\zeta}+S_{F}^{\alpha}(a)S_{F}^{\alpha}(s)^{\mu}|$ we have
\begin{align}\label{xzqa6}
  &\mathcal{L}_{F}^{\alpha,-1}\left[\frac{S_{F}^{\alpha}(s)^{\xi}}{S_{F}^{\alpha}(s)^{\zeta}+S_{F}^{\alpha}(a)S_{F}^{\alpha}(s)^{\mu}+S_{F}^{\alpha}(b)}\right]= \nonumber\\&
  S_{F}^{\alpha}(x)^{\zeta-\xi-1}\sum_{n=0}^{\infty}\sum_{k=0}^{\infty}\frac{(-S_{F}^{\alpha}(b))^n (-S_{F}^{\alpha}(a))^k}{\Gamma_{F}^{\alpha}(k(\zeta-\mu)+(n+1)\zeta-\xi)}\nonumber\\&\binom{n+k}{k}S_{F}^{\alpha}(x)^{k(\zeta-\mu)+n\zeta}.
\end{align}
\textbf{Proof:} Since we can write
\begin{align}\label{mand2}
 & \frac{S_{F}^{\alpha}(s)^{\xi}}{S_{F}^{\alpha}(s)^{\zeta}+S_{F}^{\alpha}(a)S_{F}^{\alpha}(s)^{\mu}+S_{F}^{\alpha}(b)}\nonumber\\&=
  \frac{S_{F}^{\alpha}(s)^{\xi}}{S_{F}^{\alpha}(s)^{\zeta}+S_{F}^{\alpha}(a)S_{F}^{\alpha}(s)^{\mu}}\frac{1}{1+\frac{S_{F}^{\alpha}(b)}
  {S_{F}^{\alpha}(s)^{\zeta}+
  S_{F}^{\alpha}(a)S_{F}^{\alpha}(s)^{\mu}}}\nonumber\\&=\sum_{n=0}^{\infty}\frac{S_{F}^{\alpha}(s)^\xi (-S_{F}^{\alpha}(b))^n}{S_{F}^{\alpha}(s)^{\zeta}+
  S_{F}^{\alpha}(a)S_{F}^{\alpha}(s)^{\mu}},
\end{align}
according to the  Lemma 3. the proof is complete.\\

\textbf{Some important formulas of the local fractal calculus are given below :} \cite{ab12,Golmankhaneh-107-k}:
\begin{align}
& \mathcal{L}_{F}^{\alpha}[S_{F}^{\alpha}(x)^n ]= \frac{\Gamma_{F}^{\alpha}(n+1)}{S_{F}^{\alpha}(s)^{n+1}},\nonumber\\&
\mathcal{L}_{F}^{\alpha}\left[\int_{S_{F}^{\alpha}(0)}^{S_{F}^{\alpha}(x)}f(S_{F}^{\alpha}(t))d_{F}^{\alpha}t\right]=\mathcal{L}_{F}^{\alpha}
\left[~_{0}I_{x}^{\alpha}f(S_{F}^{\alpha}(t))\right]\nonumber\\& = \frac{\mathcal{F}_{F}^{\alpha}(s)}{s},\nonumber\\&
 \mathcal{L}_{F}^{\alpha}[S_{F}^{\alpha}(x)^n f(S_{F}^{\alpha}(x))]=(-1)^{n}(D_{F}^{\alpha})^{n}\mathcal{F}_{F}^{\alpha}(s),\nonumber\\&
 \mathcal{L}_{F}^{\alpha}\left[\int_{S_{F}^{\alpha}(0)}^{S_{F}^{\alpha}(x)}f(S_{F}^{\alpha}(x)-S_{F}^{\alpha}(t))g(S_{F}^{\alpha}(t))d_{F}^{\alpha}t\right]\nonumber\\&=
\mathcal{F}_{F}^{\alpha}(S_{F}^{\alpha}(s)) \mathcal{G}_{F}^{\alpha}(S_{F}^{\alpha}(s)),
\end{align}
and
 \begin{align}\label{vcxzd}
  & \mathcal{L}_{F}^{\alpha}[(D_{F}^{\alpha})^{n}f(S_{F}^{\alpha}(x))]\nonumber\\& = (S_{F}^{\alpha}(s))^{n\alpha}
\mathcal{F}_{F}^{\alpha}(s)-(S_{F}^{\alpha}(s))^{n\alpha-1}f(S_{F}^{\alpha}(0))\nonumber\\&-
(S_{F}^{\alpha}(s))^{n\alpha-2}D_{F}^{\alpha}f(x)|_{x=S_{F}^{\alpha}(0)}-\ldots\nonumber\\&-(D_{F}^{\alpha})^{n-1}f(x)|_{x=S_{F}^{\alpha}(0)}.
 \end{align}
 \textbf{Remark 2.} If we choose  $\alpha=1$ we obtain the standard result.\\
\textbf{The important formulas of the non-local fractal calculus are as follows } \cite{Golmankhaneh-107-k}:
\begin{align}
&  _{0}\mathcal{I}_{x}^{\beta}(S_{F}^{\alpha}(x))^{\eta}=\frac{\Gamma^{\alpha}_{F}(\eta+1)}
{\Gamma^{\alpha}_{F}(\eta+\beta+1)}(S_{F}^{\alpha}(x))^{\eta+\beta},\nonumber\\&
  _{0}\mathcal{D}_{x}^{\beta}(S_{F}^{\alpha}(x))^{\eta}=\frac{\Gamma^{\alpha}_{F}(\eta+1)}
{\Gamma^{\alpha}_{F}(\eta-\beta+1)}(S_{F}^{\alpha}(x))^{\eta-\beta}.\nonumber\\&
  _{0}\mathcal{D}_{x}^{\beta}(c~\chi_{F}^{\alpha})=\frac{c}
{\Gamma^{\alpha}_{F}(1-\beta)}(S_{F}^{\alpha}(x))^{-\beta},\nonumber\\&
\mathcal{L}_{F}^{\alpha}[ _{0}\mathcal{I}_{x}^{\beta}f(x)] = \frac{\mathcal{F}_{F}^{\alpha}(S_{F}^{\alpha}(s))}
{S_{F}^{\alpha}(s)^{\beta}}.
\end{align}
where $c$ is constant.\\
\textbf{Remark 3.} If we choose $\beta=\alpha$ then we arrive at to the local fractal  derivative whose order  is equal the dimension of the fractal.\\
\section{ Comparison between the local fractal diﬀerential and non-local fractal diﬀerential \label{5-sec}}
In this section, we compare the local and non-local fractal differential equations.\\
\textbf{Example 1.} Consider  linear local fractal differential equation as
\begin{equation}\label{zswed8522w3}
   D_{F}^{\alpha} y(x)+y(x)=0,
\end{equation}
with the initial-value
\begin{equation}
y(x)|_{x=S_{F}^{\alpha}(0)}=1,~~~~
\end{equation}
\begin{figure}[H]
  \centering
  \includegraphics[scale=0.5]{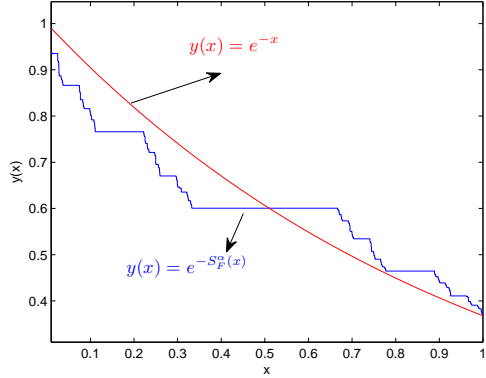}
  \caption{We plot the solution of  Eq. (\ref{zswed8522w3}).}\label{swedwqsaq67b}
\end{figure}
Hence the solution to Eq. (\ref{zswed8522w3}) is
\begin{equation}\label{zswfrddss22w3}
y(x)=e^{-S_{F}^{\alpha}(x)},
\end{equation}
where $\alpha=0.6309$ is the $\gamma$-dimension of the  triadic Cantor  set \cite{ab12,Golmankhaneh-107-k}.\\ In Figure \ref{swedwqsaq67b} we give the graph of Eq. (\ref{zswfrddss22w3}).\\
\textbf{Example 2.} Consider  linear non-local fractal differential equation as
\begin{equation}\label{xszc2586de}
      _{0}^{C}\mathcal{D}_{x}^{\beta} y(x)+y(x)=0,~~~~
\end{equation}
with the initial condition
\begin{equation}\label{cx}
y(x)|_{x=S_{F}^{\alpha}(0)}=1,~~~~D_{F}^{\alpha} y(x)|_{x=S_{F}^{\alpha}(0)}=0.
\end{equation}
In view of  Eq. (\ref{poi78})  we have
\begin{equation}\label{z}
\mathcal{L}_{F}^{\alpha}\{  _{0}^{C}\mathcal{D}_{x}^{\beta}f(x) \}=\frac{(S_{F}^{\alpha}(s))^{ \alpha}\mathcal{F}_{F}^{\alpha}(s)-1 }{S_{F}^{\alpha}(s)^{\alpha-\beta}}.
\end{equation}
\begin{figure}[H]
  \centering
  \subfloat[If we choose $\beta=0.33$ in Eq. (\ref{ujhy63cx})]{\includegraphics[width=0.5\textwidth]{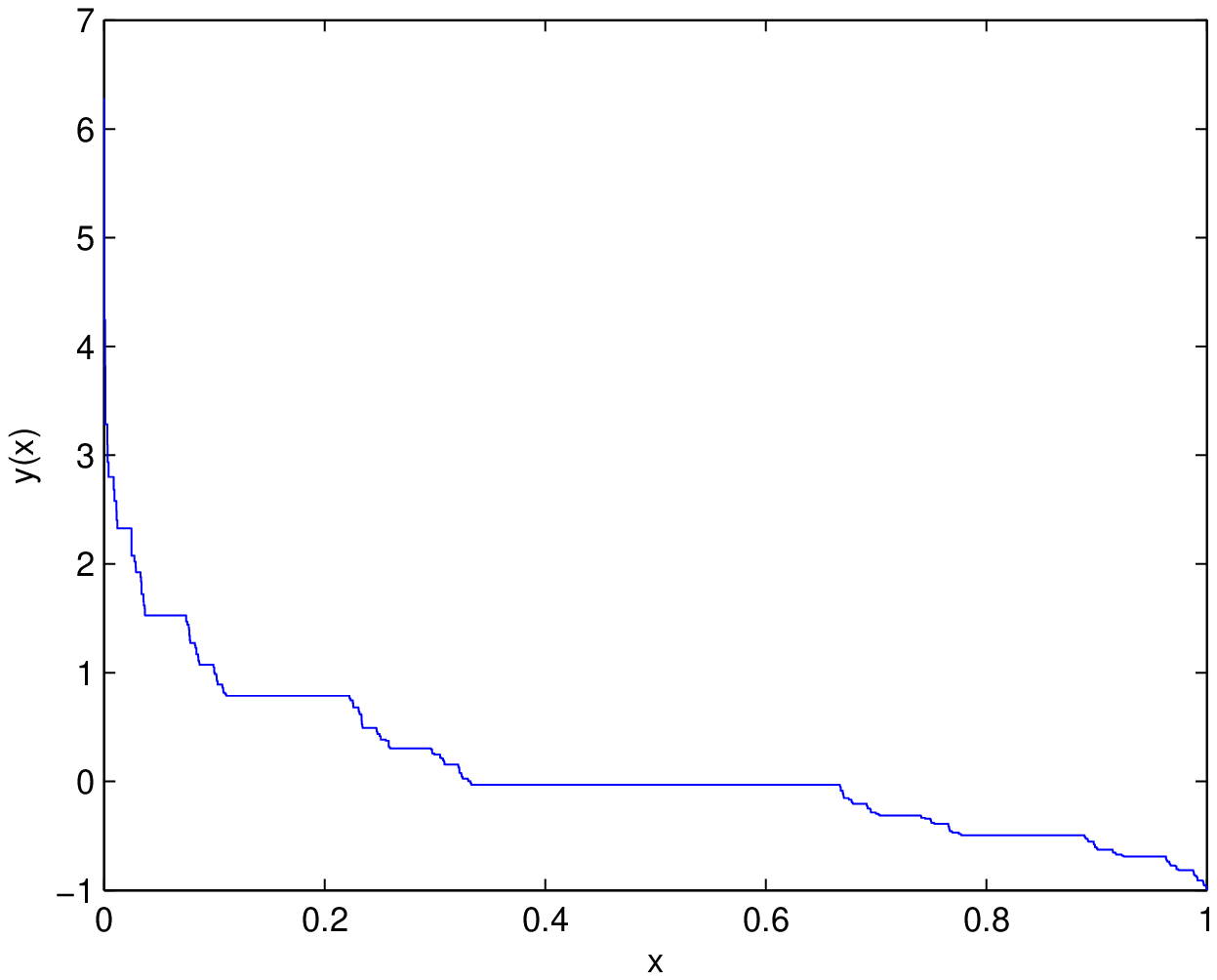}\label{fig:f1}}
  \hfill
  \subfloat[If we choose  $\beta=0.25$ in Eq. (\ref{ujhy63cx}) ]{\includegraphics[width=0.5\textwidth]{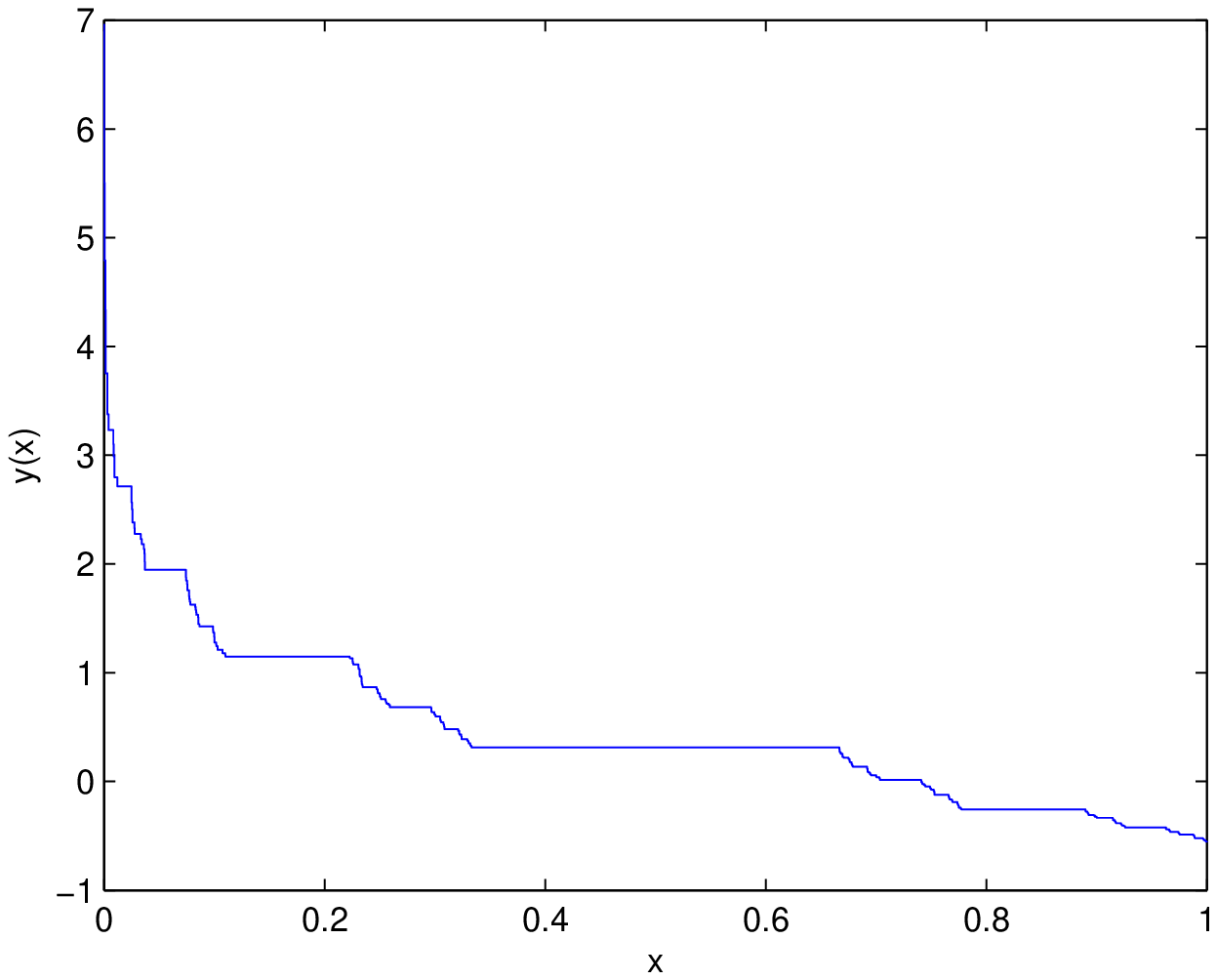}\label{fig:f2}}
  \caption{We draw the graph of  Eq. (\ref{ujhy63cx}).}\label{figh2}
\end{figure}
Applying the fractal Laplace transformation on the both sides of Eq. (\ref{xszc2586de}) and using  Eq. (\ref{poi78}) we obtain
\begin{equation}\label{b}
  \frac{(S_{F}^{\alpha}(s))^{ \alpha}\mathcal{F}_{F}^{\alpha}(s)-1 }{S_{F}^{\alpha}(s)^{\alpha-\beta}}+\mathcal{F}_{F}^{\alpha}(s)=0.
\end{equation}
It follows that
\begin{equation}\label{xzs}
 \mathcal{F}_{F}^{\alpha}(s)=\frac{S_{F}^{\alpha}(s)^{\beta-\alpha}}{1+S_{F}^{\alpha}(s)^{\beta}},
\end{equation}
using the fractal inverse Laplace transform Eq. (\ref{g895raws}) we arrive at the solution of Eq. (\ref{xszc2586de}) as follows
\begin{equation}\label{ujhy63cx}
 y(x)=S_{F}^{\alpha}(x)^{\alpha-1}E^{\alpha}_{F,\beta,\alpha}\left(-S_{F}^{\alpha}(x)^{\beta}\right).
\end{equation}
In Figure \ref{figh2} we present the graph of Eq.( \ref{ujhy63cx}).\\
\section{Application of non-local fractal differential equations \label{6-sec} }
\begin{figure}[H]
  \centering
  \includegraphics[scale=0.5]{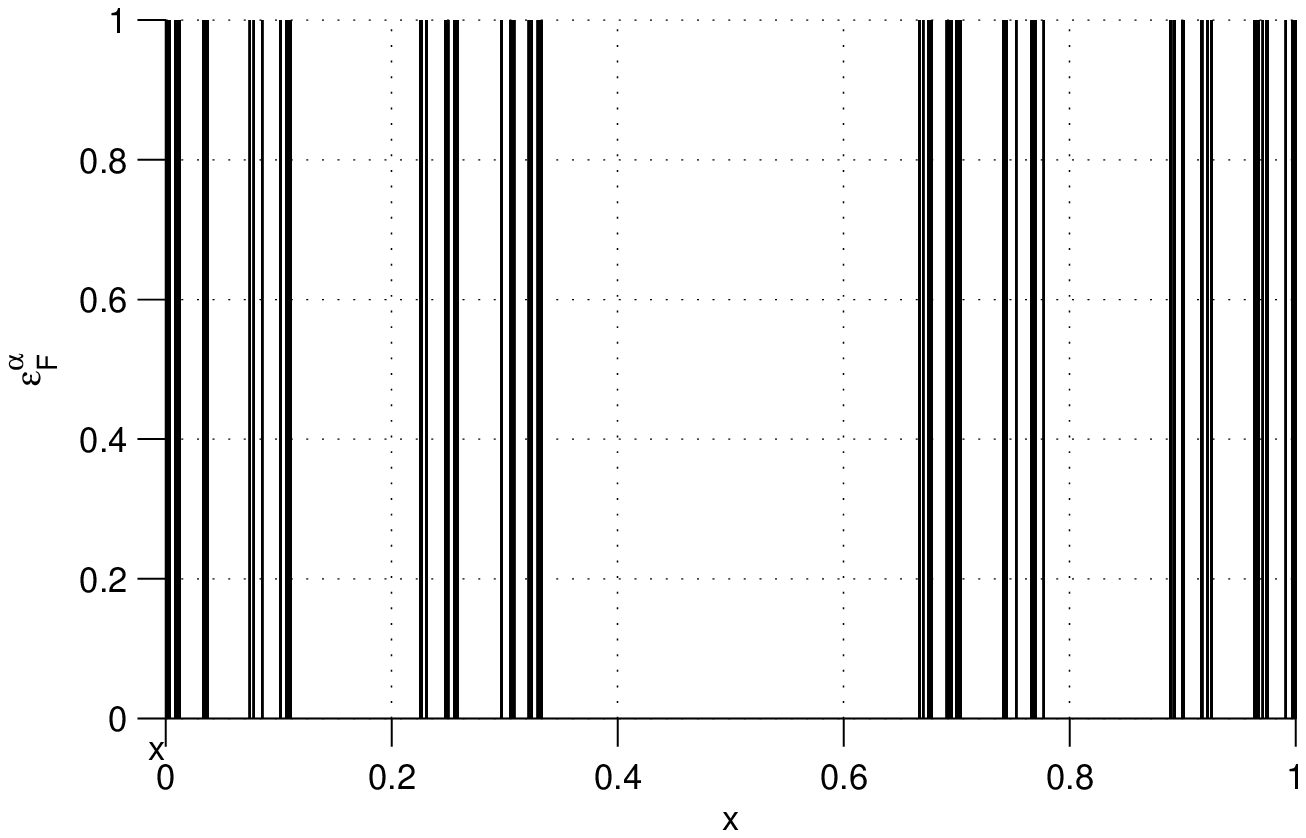}
  \caption{We sketch  $ \epsilon_{F}^{\alpha}(t)=\chi_{F}^{\alpha}$ which is characteristic function of the triadic Cantor set. }\label{751zswaq67b}
\end{figure}

 In this section we give the applications and  new models are given to the non-local fractal derivatives \cite{Golmankhaneh-107-k}.\\
\textbf{Fractal Abel's tautochrone:}  As a first example we generalized Abel's problem which is the curve  of quick descent on the fractal time-space. Using the conservation of energy in the fractal space the differential equation of the motion a particle is
\begin{equation}\label{swewaw}
  D_{F,t}^{\alpha} \mathfrak{s}_{F}^{\alpha}=\frac{d_{F}^{\alpha}\mathfrak{s}_{F}^{\alpha}}{d_{F}^{\alpha} t}=-\sqrt{2g_{F}^{\alpha}(S_{F}^{\alpha}(y)-S_{F}^{\alpha}(y_{0}))},
\end{equation}
where $\mathfrak{s}_{F}^{\alpha}$ is fractal arc length, and  $g_{F}^{\alpha}$ fractal space gravitational constant, and $y$ is the high particle from the reference of potential. As a result we have
\begin{equation}\label{vgtr}
S_{F}^{\alpha}(T)=-\frac{1}{\sqrt{2 g_{F}^{\alpha}}}\int_{S_{F}^{\alpha}(A)}^{S_{F}^{\alpha}(B)}\frac{1}{\sqrt{(S_{F}^{\alpha}(y)-S_{F}^{\alpha}(\eta))}}
d_{F}^{\alpha}\mathfrak{s}_{F}^{\alpha}.
\end{equation}
Let us consider
\begin{equation}\label{b}
  \mathfrak{s}_{F}^{\alpha}=h_{F}^{\alpha}(S_{F}^{\alpha}(\eta)),
\end{equation}
so that we have
\begin{equation}\label{vgtr}
S_{F}^{\alpha}(T)=-\frac{1}{\sqrt{2 g_{F}^{\alpha}}}\int_{S_{F}^{\alpha}(y)}^{S_{F}^{\alpha}(0)}(S_{F}^{\alpha}(y)-S_{F}^{\alpha}(\eta))^{-1/2}D^{\alpha}_{F,\eta}h_{F}^{\alpha}(\eta)
d_{F}^{\alpha}\eta.
\end{equation}
Utilizing  $D^{\alpha}_{F,\eta}h_{F}^{\alpha}(S_{F}^{\alpha}(y))=f(S_{F}^{\alpha}(y))$ we arrive at
\begin{equation}\label{vgtr}
S_{F}^{\alpha}(T)=-\frac{1}{\sqrt{2 g_{F}^{\alpha}}}\int_{S_{F}^{\alpha}(y)}^{S_{F}^{\alpha}(0)}(S_{F}^{\alpha}(y)-S_{F}^{\alpha}(\eta))^{-1/2}f(S_{F}^{\alpha}(y))
d_{F}^{\alpha}\eta.
\end{equation}
It follows
\begin{equation}\label{bvlkcx}
\frac{\sqrt{2 g_{F}^{\alpha}}}{\Gamma(\frac{1}{2})}  S_{F}^{\alpha}(T)=~_{0}\mathcal{D}_{y}^{1/2} f(y).
\end{equation}
 The solution of Eq.(\ref{bvlkcx}) is called the  fractal cycloid.\\
 \textbf{Fractal models for the viscoelasticity:} We generalize the viscoelasticity models to the fractal mediums in the case of ideal solids and ideal liquids.~Namely, the fractal  ideal solids  describe by
 \begin{equation}\label{xc987f}
   \sigma_{F}^{\alpha}(t)=E_{F}^{\alpha}\epsilon_{F}^{\alpha}(t),
 \end{equation}
 which is called  Hooke's Law of fractal elasticity. Where $\sigma_{F}^{\alpha}$ is fractal stress, $\epsilon_{F}^{\alpha}$ is fractal strain which occurs under the applied stress and $E_{F}^{\alpha}$ is the elastic modulus of the fractal material.\\
 The fractal ideal fluid can model and describe  by  Newton's Law of fractal viscosity  as follows
 \begin{equation}\label{963u}
  \sigma_{F}^{\alpha}(t)=\lambda_{F}^{\alpha}~D_{F}^{\alpha} \epsilon_{F}^{\alpha}(t),
 \end{equation}
 where $\lambda_{F}^{\alpha}$ is the viscosity of the fractal material. But in the nature we have real martials which have properties between the ideal solids and ideal liquids. It is clear that  in the Hooke's Law of fractal elasticity Eq. (\ref{xc987f}) fractal stress is proportional to the $0$-order derivative of the fractal strain and in the Newton's Law of fractal viscosity the stress is proportional to the $\alpha$-order derivative of the fractal strain. Therefore, more general  model is
 \begin{equation}\label{wqer39}
   \sigma_{F}^{\alpha}(t)=E_{F}^{\alpha} (\chi_{F}^{\alpha})^{\beta}~ _{0}\mathcal{D}_{x}^{\beta} \epsilon_{F}^{\alpha}(t),~~~~~~~~\chi_{F}^{\alpha}=\frac{\lambda_{F}^{\alpha}}{E_{F}^{\alpha}},
 \end{equation}
which is called fractal Blair's model.~Here, we suggest the fractional non-local order fractal derivative  $\beta$ as an index of memory. Namely, if we choose $\beta=0$ in the process is nothing forgotten  and the case of $\beta=\alpha$  the process is memoryless. Hence if we choose  $0<\beta<\alpha$ it shows the processes with memory on the fractals.\\ If we choose
\begin{equation}\label{cd}
  \epsilon_{F}^{\alpha}(t)=\chi_{F}^{\alpha},
\end{equation}
where $\chi_{F}^{\alpha}$ is characteristic function of the triadic Cantor set. In Figure \ref{751zswaq67b} we plot the $\epsilon_{F}^{\alpha}(t)$.
Utilizing Eq. (\ref{wqer39}) we obtain the  fractal stress as follows
\begin{equation}\label{olp}
  \sigma_{F}^{\alpha}(t)=E_{F}^{\alpha} (\chi_{F}^{\alpha})^{\beta}\frac{1}
{\Gamma^{\alpha}_{F}(1-\beta)}(S_{F}^{\alpha}(t))^{-\beta}.
\end{equation}
\begin{figure}[H]
  \centering
  \includegraphics[scale=0.5]{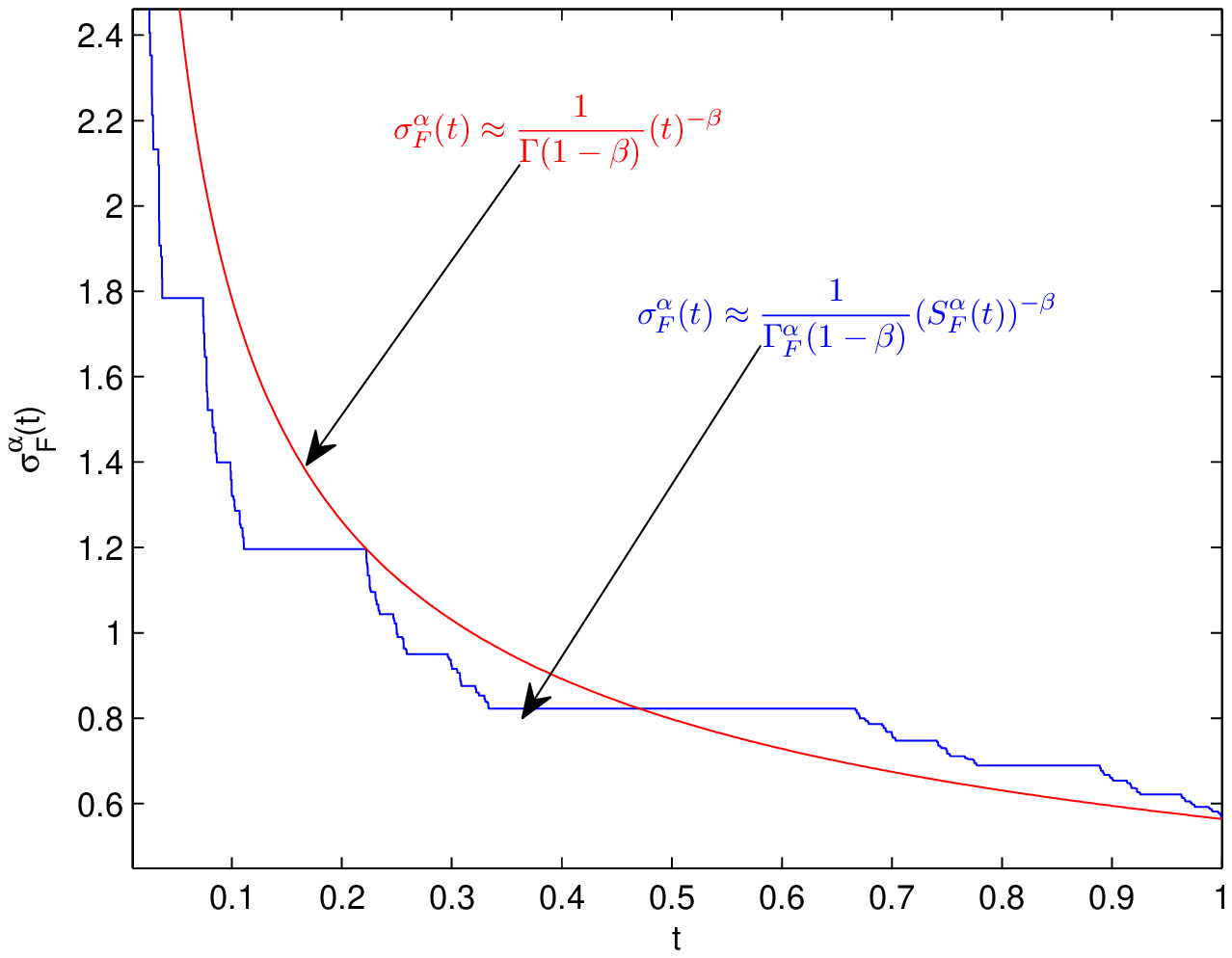}
  \caption{We sketch $\sigma_{F}^{\alpha}(t)$ for the fractal stress substituting  $\beta=0.5$ }\label{751zkilujswaq67b}
\end{figure}
In Figure \ref{751zkilujswaq67b} we show the graph of $\sigma_{F}^{\alpha}(t)$ fractal stress.\\
\textbf{Remark 4.}  If we choose $\beta=0$ and $\beta=\alpha$ in Eq. (\ref{wqer39}) we will  have the fractal stress and  the fractal strain relations for the cases of fractal ideal solids and  the fractal ideal fluids, respectively.
\section{ Conclusion \label{7-sec}}
In this paper we generalized the fractal calculus involving the non-local derivatives. The scaling properties of the local and non-local derivatives are studied  because they are important in   physical applications. Using an illustrative example we compared the local and non-local linear fractal differential equations. We also suggested  some applications for the new non-local fractal differential equations.
%\begin{thebibliography}{99}
%\bibitem{ref1} ...
%\end{thebibliography}

\end{document}